\documentclass[]{amsart}

\usepackage{amssymb,amsfonts, amsmath, mathrsfs}

\title[fixed points theorems for sum of operators]{ A Topological and
Geometric Approach to Fixed Points Results for Sum of Operators and
Applications }

\author{Cleon S. Barroso\and Eduardo V. Teixeira}

\thanks{The first author is grateful for the financial support by Instituto
do Mil\^enio - CNPq}
\thanks{The second author acknowledges support from CNPq}

\address{Departamento de Matem\'atica, Universidade Federal do Cear\'a,
Campus do Pici, Bl. 914, 60455-760,
Fortaleza, CE, Brazil.}

\email{cleonbar@mat.ufc.br}

\address{Department of Mathematics, University of Texas at Austin, RLM
9.136, Austin, Texas 78712-1082.}
\email{teixeira@math.utexas.edu}
\newtheorem{theorem}{Theorem}[section]

\newtheorem{lemma}[theorem]{Lemma}

\newtheorem{proposition}[theorem]{Proposition}

\newtheorem{corollary}[theorem]{Corollary}

\theoremstyle{definition}

\newtheorem{definition}[theorem]{Definition}

\newtheorem{example}[theorem]{Example}

\theoremstyle{remark}

\newtheorem{remark}[theorem]{Remark}

\numberwithin{equation}{section}

\begin{document}


\begin{abstract}
In this paper we establish a fixed point result of Krasnoselskii
type for the sum $A+B$, where $A$ and $B$ are continuous maps
acting on locally convex spaces. Our results extend previous ones.
We apply such results to obtain strong solutions for some
quasi-linear elliptic equations with lack of compactness. We also
provide an application to the existence and regularity theory of
solutions to a nonlinear integral equation modeled in a Banach
space. In the last section we develop a sequentially weak
continuity result for a class of operators acting on vector-valued
Lebesgue spaces. Such a result is used together with a geometric
condition as the main tool to provide an existence theory for
nonlinear integral equations in $L\sp p(E)$.
\end{abstract}

\maketitle

\section{Introduction}

Many problems arising from the most diverse areas of natural
science, when modeled under the mathematical point of view,
involve the study of solutions of nonlinear equations of the form
\begin{equation}\label{eqn:0}
Au+Bu=u,\quad u\in M,
\end{equation}
where $M$ is a closed and convex subset of a Banach space $X$, see
for example \cite{1,CSBarroso,2,Dhage1,Dhage2,Dhage3}. In
particular, many problems in integral equations can be formulated
in terms of (\ref{eqn:0}). Krasnoselskii's fixed point Theorem
appeared as a prototype for solving equations of the type
(\ref{eqn:0}), where $A$ is a continuous and compact operator and
$B$ is, in some sense, a contraction mapping. Motivated by the
observation that the inversion of a perturbed differential
operator could yield a sum of a contraction and a compact
operator, Krasnoselskii proved that the sum $A+B$ has a fixed
point in $M$, if: {\bf (i)} $A$ is continuous and compact, {\bf
(ii)} $B$ is a strict contraction and {\bf (iii)} $Ax+By\in M$ for
every $x,y\in M$. Since then a wide class of problems, for
instance in integral equations and stability theory, have been
contemplated by the Krasnoselskii fixed point approach. However,
in several applications, the verification of {\bf (iii)} is, in
general, quite hard or even impossible to be done. As a tentative
approach to grapple with such a difficulty, many interesting works
have appeared in the direction of relaxing hypothesis {\bf (iii)}.
\par In a recent paper \cite{1}, Burton proposes the following improvement for {\bf
(iii)}: $(\,\mbox{If } u=Bu+Av\, \mbox{ with }v\in M, \mbox{ then
} u\in M). $ Subsequently, in \cite{3}, the following new
asymptotic requirement was introduced:
$$
\big(\,\mbox{If } \lambda\in (0,1) \mbox{ and } u=\lambda Bu+Av\,\mbox{ for
some } v\in M,\mbox{ then } u\in
M\big),
$$
In this paper we explore another kind of generalization. Indeed we
study suitable modifications on conditions  {\bf (i)} and  {\bf
(ii)} as well. Notice that condition  {\bf (i)} involves
continuity and compactness. Based on the well known fact that
infinite dimensional Banach spaces are not locally compact, we
suggest a locally convex topology approach to equation
(\ref{eqn:0}). The interpretation of some practical equation of
the form (\ref{eqn:0}) may face the problem that the operators
involved are not even continuous. The freedom of choosing a more
general notion of topology might remedy this difficulty. We should
mention that others authors have already studied equation
(\ref{eqn:0}) in
locally convex spaces \cite{2,V}. \\
\indent Condition {\bf (ii)} is also, in same sense, quite
restrictive. Indeed this condition implies norm-continuity. In
this paper we shall suggest a much more general condition than
strict contraction. All the generality of our results will be used
in the applications. For instance, in section $6$, when we shall
be interested in proving optimal regularity of solutions to a
nonlinear integral equation, we will apply our fixed point results
on the space $W^{1,\infty}(I,E)$. The operators in question will
be neither norm continuous nor weakly continuous. Thus the
classical assumptions on fixed point results of Krasnoselskii type
does not hold.

\par Our paper is organized as follows: In Section 2  we reformulate the
Krasnoselskii fixed point Theorem for the locally convex setting.
This new version of the Krasnoselskii fixed point Theorem we
provide here generalizes, among others, the results in \cite{3}.
The main abstract results are introduced in this section: Theorem
\ref{trm:I}, Corollary \ref{cor:1} and Theorem \ref{trm:1}.

\par In Section 3, we apply our fixed point theory to the solvability of one
parameter operator equations of the form
$$
    Au + \lambda Bu = u,~\lambda \ge 0.
$$
In this section we shall restrict ourself to reflexive Banach
spaces endowed with the weak topology. In the next section we
exemplify the power of results established by studying an elliptic
equation with lack of compactness. In particular we solve a
quasi-linear elliptic equation with Sobolev critical exponent. \\
\indent A nonlinear integral equation is studied in section 5. We
provide an existence principle for the following nonlinear
integral equation:
\begin{equation}\label{eqn:0.1}
    u(t)=f(u(t))+\int_0\sp tg(s,u(s))ds, \quad u\in C(I,E),
\end{equation}
where $E$ is a reflexive Banach space and $I=[0,T]$. In section 6
we explore the optimal regularity of solutions of equation
(\ref{eqn:0.1}). The approach in this section explores all the
generality of the results established in section 2. As part of our
strategy to get Lipschitz regularity for solutions of equation
(\ref{eqn:0.1}), we suggest a new locally convex topology for the
space $W^{1,\infty}(I,E)$. In such topology the ball of this space
is compact. Furthermore, we can proof continuity of all the
operators involved in equation (\ref{eqn:0.1}). Such operators are not norm-continuous
 though. \\
\indent Finally, in Section 7, we explore a new geometric idea of
finding fixed point for sum of operators. Such an approach is
motivated by a sort of strong triangular inequality for uniformly
convex spaces. We apply this idea directly in the study of the
following challenging variant of (\ref{eqn:0.1}):
\begin{equation}\label{eqn:0.2}
    u(t)=f(t,u(t))+\Phi \left (t, \int_0\sp t k(t,s)u(s)ds \right ),\quad
u\in L\sp p(I,E),
\end{equation}
where $E$ is a uniformly convex space, $1<p<\infty$ and $I=[0,T]$.
As usual, the existence of solutions to (\ref{eqn:0.2}) reduces to
search fixed points for the operator $A+B$, where
\begin{eqnarray*}
&&Au(t)=f(t,0)+\Phi\left (t,\int_0\sp t k(t,s)u(s)ds\right ),\\
&&Bu(t)=f(t,u(t))-f(t,0).
\end{eqnarray*}
A geometric condition is used in order to assure that
$(A+B)(B\sb{L\sp p(E)}(\overline{R}))\subseteq B\sb{L\sp
p(E)}(\overline{R})$, for some $\overline{R}>0$. Our final step
toward the solution of  equation (\ref{eqn:0.2}) is a sequentially
weak continuity result (Lemma \ref{weak continuity}) which
guarantees the weak sequential continuity for the operator acting
on vector-valued Lebesgue spaces.


\section{Abstract Fixed Point Theorems for Sum of Operators}
The notation and terminology used in this paper are standard. For
convenience of the reader we recall some basic facts. Let
$(X,\mathcal{T})$ be a Hausdorff locally convex topological vector
space. The symbol $\mathcal{T}$ stands for the family of seminorms
that generates the locally convex topology on $X$. In the sequel
we will use the following well-known Schauder-Tychonoff fixed
point Theorem for locally convex spaces.

\begin{theorem}\label{Schauder} Let $M$ be a closed convex subset of a
Hausdorff
locally convex space $X$ and let $T:M\to M$ be a continuous mapping such
that $T(M)$
is relatively compact. Then $T$ has a fixed point in $M$.
\end{theorem}

Our first version of the Krasnoselskii fixed point Theorem is as
follows.
\begin{theorem}\label{trm:I} Let $M$ be a closed  and convex subset of
a Hausdorff locally convex space $X$ and $A,B:M\to X$ be continuous
operators such
that
\begin{itemize}
\item[(a)] $A(M)$ is relatively compact;

\item[(b)] $(I-B)$ is continuously invertible and $A(M) \subseteq (I-B)(M)$;

\item[(c)] If $u = B(u) +A(v)$ for some $v \in M$ then $u \in M$.
\end{itemize}
Then $A+B$ has a fixed point in $M$.
\end{theorem}

\begin{proof} Let us define $T\colon M \to M$ by
$$
    T(u) := (I-B)^{-1} \circ A(u).
$$
$T$ is well defined by item $(a)$. The fact that $T$ maps $M$ into $M$
follows by
condition $(c)$. Furthermore $T$ is a continuous map and $T(M)$ is
relatively
compact. Applying Theorem \ref{Schauder} to the operator $T$ we conclude
there exists
a $u \in M$ such that $T(u) = u$. Finally notice that a fixed point to $T$
is
actually a fixed point to $A+B$.
\end{proof}

\begin{corollary} \label{cor:1} Let $M$ be a compact convex subset of
a Hausdorff locally convex space $X$ and $A,B:M\to X$ be continuous
operators such
that
\begin{itemize}
\item[(a)] There exists a sequence $\lambda_n \to 1$, such that
$(I-\lambda_n B)$ is
injective and $A(M) \subseteq (I - \lambda_n B)(M),~\forall n$;

\item[(b)] For all $n\ge 1$, if  $u = \lambda_n B(u) +A(v)$ for some $v \in
M$ then
$u \in M$.
\end{itemize}
Then $A+B$ has a fixed point in $M$.
\end{corollary}
\begin{proof} We first notice that once $M$ is compact and $B$ is
continuous, we automatically have the continuity of $(I-\lambda_n
B)^{-1}$ and hence of the operators $T_n := (I-\lambda_n B)^{-1}
\circ A$ for each $n \ge1$. Indeed, let $\xi_\gamma$ be a net in
$M$ converging to $\xi \in M$. Let us denote by $\psi_\gamma =
(I-\lambda_n B)^{-1} \circ A(\xi_\gamma)$. Since $M$ is compact,
up to a subnet we might assume $\psi_\gamma \to \psi$, thus, using
the continuity of $A$ and $B$ we have
$$
    A(\xi_\gamma) = (I - \lambda_n B)(\psi_\gamma) \to (I - \lambda_n
    B)(\psi).
$$
Thus $\psi = (I - \lambda_n B)^{-1} \circ A (\xi)$. We now can
apply, for each $n \ge 1$, Theorem \ref{trm:I} to $A$ and
$\lambda_n B$ in order to get a fixed point $u_n$ for $A +
\lambda_n B$, i.e, there exists a $u_n$ such that
\begin{equation} \label{FPn}
    u_n = \lambda_n B(u_n) + A(u_n).
\end{equation}
Finally, up to a subnet we may suppose $u_n \to u$ in $M$. Passing the limit
in
(\ref{FPn}) we conclude the proof of Corollary \ref{cor:1}.
\end{proof}
Another variant of Corollary \ref{cor:1}, where the compactness of
$M$ is substituted by the relative compactness of $A(M)$ is the
following.
\begin{corollary} \label{cor:2} Let $M$ be a convex and closed
subset of a Hausdorff locally convex space $X$ and $A,B \colon M
\to X$ be continuous operators such that
\begin{itemize}
\item[(a)] $(I- B)$ is injective and $A(M) \subseteq (I - \lambda_n
B)(M),~\forall n$;

\item[(b)] $A(M)$ is relatively compact.

\item[(c)] If  $u = \lambda_n B(u) +A(v)$ for some $v \in
M$ then $u \in M$.
\end{itemize}
Then $A+B$ has a fixed point in $M$.
\end{corollary}

\begin{remark} It is worthwhile to highlight that the condition
of $A(M) \subseteq (I - \lambda B)(M)$ in Theorem \ref{trm:I},
Corollary \ref{cor:1} and Corollary \ref{cor:2} can be relaxed by
$A(M) \subseteq (I - \lambda B)(\mathcal{D}(B))$, where
$\mathcal{D}(B)$ stands for the biggest domain for which $B$ is
continuous on it.
\end{remark}

 A simple way of checking the invertibility of $(I - B)$
and the condition of $A(M) \subseteq (I-B)(M)$ is to ask that $B$
is a contraction is the following sense.
\begin{definition} Let $B \colon M \to X$ be an operator defined
in a subset of a locally convex space $X$. Let $\mathcal{T}$ be a family of
seminorms
that define the topology in $X$. We say $B$ is a $\mathcal{T}$-contraction
if for
each $\rho \in \mathcal{T}$ there exists a $\lambda_\rho \in (0,1)$ such
that
$$
    \rho( B(u) - B(v)) \le \lambda_\rho \rho(u-v).
$$
\end{definition}
Indeed when $X$ is complete, it is a simple adaptation of the
original proof of the Banach fixed point Theorem the fact that
every map $B \colon X \to X$ has a unique
fixed point. For more details see \cite{2}. \\
\indent Let now $(X,\|\cdot\|)$ be a Banach space and let
$\mathcal{T}$ be the family of seminorms $\{\rho\sb f(x)=|\langle
f,x\rangle|\, :\, f\in X\sp * \mbox{ and } \|f\|\sb{X\sp *}\leq 1
\}$. The topology generated by $\mathcal{T}$ is called the weak
topology. In practical situations, often one faces the problem of
solving equations of the type (\ref{eqn:0}) in the weak topology
setting. One of the advantages of this special locally convex
topology is the fact that if a set $M$ is weakly compact, then
every sequentially weakly continuous map $T:M\to X$, i.e. an
operator which maps weakly convergent sequences into weakly
convergent sequences, is weakly continuous. This is an immediately
consequence of the fact that weak sequential compactness is
equivalent to weak compactness (Eberlein-$\check{S}$mulian's
Theorem). From this observation, the following version of Schauder
fixed point Theorem holds \cite{5}.
\begin{theorem}\label{Schauder II} Let $M$ be a convex weakly compact subset
of a Banach space $X$. Then every sequentially weakly continuous
operator $T$ self-mapping $M$ has a fixed point.
\end{theorem}

Consequently we have the following result which will be used in
this form in section 5. Such a result encloses an improvement
which will turn out to be crucial for establishing the existence
principle for the nonlinear integral equation studied in section
5. Before stating this result we need a definition.
\begin{definition} Let $M$ be a closed and convex subspace of a
Hausdorff locally convex space $X$ and $A,B \colon M \to X$
continuous operators. We will denote by $\mathcal{F} =
\mathcal{F}(M,A,B)$ the following set
$$
    \mathcal{F} := \{ u \in X : u = B(u) + A(v) \textrm{ for some
    } v \in M \}.
$$
\end{definition}

\begin{theorem}\label{trm:1} Let $M$ be a closed, convex subset of a Banach
space $X$.
Assume that $A,B:M\to X$ satisfies:
\begin{itemize}
\item[(a)] $A$ is sequentially weakly continuous;

\item[(b)] $B$ is $\lambda$-contraction;

\item[(c)] If $u=Bu+Av$, for some $v\in M$, then $u\in M$;

\item[(d)] If $\{u\sb n\}$ is a sequence in $\mathcal{F}$ such that $u\sb
n\rightharpoonup u$, for some $u\in M$, then $Bu\sb n\rightharpoonup Bu$;

\item[(e)] The set $\mathcal{F}$ is relatively weakly compact.
\end{itemize}
Then $A+B$ has a fixed point in $M$.
\end{theorem}

\begin{proof} Fix a point $u\in M$ and let $Tu$ be the unique point in $X$
such that $Tu=B\cdot Tu+Au$. By (c), we have $Tu\in M$. So that
the mapping $T:M\to M$ given by $u\mapsto Tu$ is well-defined.
Notice that $Tu=(I-B)\sp{-1}Au$, for all $u\in M$. In addition, we
observe that $T(M)\subset \mathcal{F}$. We claim now that $T$ is
sequentially weakly continuous in $M$. Indeed, let $\{u\sb n\}$ be
a sequence in $M$ such that $u\sb n\rightharpoonup u$ in $M$.
Since $\{Tu\sb n\}\subset\mathcal{F}$, the assumption (e)
guarantees, up to a subsequence,  that $Tu\sb n\rightharpoonup v$,
for some $v\in M$. By (d), we have $B\cdot Tu\sb n\rightharpoonup
Bv$. Also, from (a) it follows that $Au\sb n\rightharpoonup Au$
and hence the equality $Tu\sb n=B\cdot Tu\sb n+Au\sb n$ give us
$v=Bv+Au$. By uniqueness, we conclude that $v=Tu$. This proves the
claim. Take now the subset
$C=\overline{\rm{co}}(\mathcal{F})\subset M$. Krein-\v{S}mulian
Theorem implies that $C$ is a weakly compact set. Furthermore, it
is easy to see that $T(C)\subset C$. Applying Theorem
\ref{Schauder II}, we find a fixed point $u\in C$ for $T$.
Consequently, this proves Theorem \ref{trm:1}.
\end{proof}

Let us now state some other consequences of Theorem \ref{trm:1}.
The first one is the following result for reflexive Banach spaces,
where closed, convex and bounded sets are weakly compacts.

\begin{corollary}\label{crl:2} Assume the conditions $(a)$-$(d)$ of Theorem
\ref{trm:1} for $A$ and $B$. If $M$ is a closed, convex and bounded
subset of a reflexive Banach space, then $A+B$ has a fixed point in $M$.
\end{corollary}

Next we consider the case when $B$ is a nonexpansive mapping (or
$1$-Lipschitz mapping) on $X$, that is a mapping satisfying
$\|Bu-Bv\|\leq \|u-v\|$, for all $u,v\in X$.

\begin{corollary}\label{crl:3} Let $M$ be a convex and weakly compact subset
of a Banach
space $X$ and let $A,B:M\to X$ be sequentially weakly continuous operators
such that
\begin{itemize}
\item[(a)] $B$ is nonexpansive;

\item[(b)] If $\lambda\in (0,1)$ and $u=\lambda Bu+Av$ with $v\in M$, then
$u\in M$;

\end{itemize}
Then $A+B$ has a fixed point in $M$.
\end{corollary}

The above versions of Krasnoselskii fixed point Theorem generalize
among others, \cite{2}, \cite{3}.


\section{Local Versions of Krasnoselskii Fixed Point Theorem \\to One Parameter Operator
Equations}

Let $X$ be a Banach space. The main goal of this section is to present some
existence
results for the following nonlinear equation operator on Banach spaces.
\begin{equation}\label{eqn:opeqn}
Au+\lambda Bu=u,
\end{equation}
where $A,B:X\to X$ and $\lambda\geq 0$. In order to do so, we
establish some local versions of above abstract results. We recall
that a mapping $T:X\to X$ is called Lipschitzs if $\|Tu-Tv\|\leq
\|T\|\sb{Lip}\|u-v\|$, for all $u,v\in X$, where $\|T\|\sb{Lip}$
denotes the Lipschitzs constant of $T$.

In the sequel we need of the following definition.

\begin{definition}\label{def:1} A mapping $T:X\to X$ is said to be expanding
if $\|u\|\leq
\|u-\lambda Tu\|$, holds for any $\lambda >0$ and all $u\in X$.
\end{definition}

Our first result concerning about equation (\ref{eqn:opeqn}) is as
follows.

\begin{theorem}\label{trm:3.2} Let $X$ be a reflexive Banach space. Assume
that
$A,B:X\to X$ are sequentially weakly continuous maps satisfying:
\begin{itemize}
\item[\rm{(i)}] $A(B\sb R)\subset B\sb R$, for some $R>0$; \item[\rm{(ii)}]
$B$ is
Lipschitzs and expanding on $X$.
\end{itemize}
Then (\ref{eqn:opeqn}) is solvable for all $\lambda\geq 0$.
\end{theorem}

\paragraph{Proof} Given $\lambda \geq 0$, we define
$\mathcal{A},\mathcal{B}:X\to X$
by
$$
\mathcal{A}(u)=\dfrac{Au+\lambda \|B\|\sb{Lip}\cdot
u}{1+\lambda\|B\|\sb{Lip}}\quad\mbox{ and }\quad
\mathcal{B}(u)=\dfrac{\lambda
Bu}{1+\lambda\|B\|\sb{Lip}}.
$$
Then, $\mathcal{A}$, $\mathcal{B}$ are sequentially weakly
continuous maps, $\mathcal{A}$ maps $B\sb R$ into itself and
$\mathcal{B}$ is a $\frac{\lambda \|B\|\sb{Lip}}{1+\lambda
\|B\|\sb{Lip}}$-contraction. Now, since $X$ is reflexive, it
follows that $B_R$ is weakly compact. Thus, by Corollary
\ref{cor:1}, $\mathcal{A}+\mathcal{B}$ has a fixed point $u\in
B_R$. Now, one easily verifies that $Au+\lambda Bu=u$. This
completes the proof.

Now we can state and prove our second result related o the
solvability of equation (\ref{eqn:opeqn}).

\begin{theorem}\label{trm:3.1} Let $X$ be a reflexive Banach space and
$B:X\to X$ a Lipschitz mapping which is sequentially weakly
continuous. Suppose that for any $\mu>0$, $A\sb\mu :X\to X$ is a
sequentially weakly continuous mapping such that
\begin{itemize}
\item[\rm{(i)}] $\|A\sb\mu u\|\leq \mu \|u\|\sp p+a\|u\|\sp q+b$,
\end{itemize}
where $p>1$, $0<q<1$ and $a,b>0$. Then there exists $\mu\sp *>0$
such that for any $\mu\in (0,\mu\sp *)$ and any $0\leq \lambda\leq
1/\|B\|_{\textrm{Lip}}$, the sum $A\sb\mu +\lambda B$ has a fixed
point.
\end{theorem}

\begin{proof} We might without loss of generality suppose $B(0) = 0$.
Fix $0<\lambda<\frac{1}{\|B\|}$. Consider the ball $M=B\sb
R(0)$ of $X$, where $R>0$ is such that
$$
\frac{a}{R\sp{1-q}}+\frac{b}{R}\leq (1-\lambda \|B\|).
$$
Now taking $\mu\sp *>0$ such that
\begin{equation}\label{eqn:mu}
\mu\sp * R\sp p +aR\sp q+\lambda \|B\|R+b\leq R,
\end{equation}
we conclude that for any $\mu\in (0,\mu\sp *)$, the sum $A\sb\mu +\lambda B$
maps $M$
into itself. Since $M$ is weakly compact, we can apply Theorem \ref{Schauder
II} to
get a fixed point to $A\sb\mu +\lambda B$. This complete the proof.
\end{proof}

In our next result we consider the case where
$B\equiv\mbox{constant}$. The proof we shall present follows the
arguments in \cite{Reichbach}.

\begin{theorem}\label{trm:3.4} Let $X$ be a Banach space and $A:X\to X$ a
compact operator such that
\begin{itemize}
\item[\rm{(i)}] $\|Au\|\leq a\|u\|\sp p$,
\end{itemize}
where $a>0$ and $p>1$. Then, there exists $R>0$ such that for any $h\in
B_R$,
$(\ref{eqn:opeqn})$ has a solution with $B\equiv h$.
\end{theorem}

\paragraph{\it Proof} For each $r>0$ let $\delta\sb r$ be the number given
by
$$
\delta\sb r=\sup\sb{\|x\|\leq r}\|Ax\|.
$$
From the assumption $\rm{(i)}$, we can choose $\sigma>0$ such that
\begin{equation}\label{eqn:3}
\inf\sb{0<r<\sigma}\dfrac{\delta\sb r}{r}<1.
\end{equation}
By (\ref{eqn:3}), there exists $r>0$ such that if $\|x\|\leq r$ we have
$$
\|Ax\|\leq \delta\sb r<r.
$$
We now define the map $T:X\to X$ by $Tx=Ax+h$. Thus, taking
$0<R<r-\delta\sb r$ and any $h\in B\sb R$, we have that
\begin{equation}\label{eqn:4}
\|Tx\|\leq \|Ax\|+\|h\|\leq r,
\end{equation}
for all $x\in B_r$. This tells us that $A$ maps $B_r$ into itself.
By the Schauder fixed point Theorem, there exists $u\in B\sb r$
such that $u=Au+h$.

\begin{remark} Other local versions of Krasnoselskii fixed point Theorem
type can be found in \cite{Dhage1}.
\end{remark}

\section{An Elliptic Equation with Lack of Compactness}

The main purpose of this section is to apply the above results in order to
study the
existence of strong solutions for a class of nonlinear elliptic problems of
the form
\begin{equation}\label{eqn:5}
-\Delta u+\lambda u=f(x,u,\mu) \mbox{ in } \Omega,\quad u=0 \mbox{ on }
\partial\Omega,
\end{equation}
where $\Omega\subset\mathbb{R}\sp N$ $(N\geq 3)$ is a bounded
domain with $C\sp{1,1}$-boundary $\partial\Omega$, $\lambda$ is a
real number and
$f:\Omega\times\mathbb{R}\times\mathbb{R}\sp{+}\to\mathbb{R}$. As
usual, by a strong solution of (\ref{eqn:5}) we mean a function
$u\in W\sp{2,2}(\Omega)\cap W\sb 0\sp{1,2}(\Omega)$ satisfying
(\ref{eqn:5}) in the sense almost everywhere. Here,
$W\sp{2,2}(\Omega)$ is the usual Sobolev space and $W\sb
0\sp{1,2}(\Omega)$ is the closure of $C\sb 0\sp{\infty}(\Omega)$
in the norm $\|u\|\sb{1,2,\Omega}$. The basic idea is to reduce
(\ref{eqn:5}) to a fixed point problem in $L\sp 2(\Omega)$ of the
form
\begin{equation}\label{eqn:fpprb}
u=N\sb{f_\mu}\circ L\sp{-1}(u)-\lambda L\sp{-1}(u),
\end{equation}
where $N\sb{f_\mu}$ is the Nemytskii operator associated to $f$ and
$L\sp{-1}$ is the
inverse of $L=-\Delta$. Thus, a solution of (\ref{eqn:fpprb}) will be a
strong
solution to (\ref{eqn:5}).

For convenience of the reader, before stating our existence result
to (\ref{eqn:5}) we recall some basic facts which will used later.
It is well-known that the mapping $u\mapsto Lu$ is one-to-one from
$W\sp{2,2}(\Omega)\cap W\sb 0\sp{1,2}(\Omega)$ onto $L\sp
2(\Omega)$. Moreover, there exists $C>0$ such that
\begin{equation}\label{eqn:6}
\|u\|\sb{2,2}\leq C\|Lu\|\sb 2,
\end{equation}
for all $u\in E$, see \cite{GT}.
\par
In what follows we consider the following assumption:
$$
(H\sb p)\left\{
\begin{array}{ccc}
&&p>1/2, \mbox{ if } N=3,4, \\
&&1<p\leq N/(N-4), \mbox{ if } N>4.
\end{array}\right.
$$
The next basic $L\sp p$-estimate will be used in the sequel.
\begin{proposition} Assume $(H\sb p)$. Then for any $w\in E$
\begin{equation}\label{eqn:7}
\|w\|\sb{2p}\leq\gamma \|Lw\|\sb 2,
\end{equation}
where $\gamma$ is a constant depending only on $p,C$ and $\Omega$.
\end{proposition}

\begin{proof} This is a consequence of the Sobolev inequalities. Assume
$1/2<p<\infty$. If
$N=3$, then it follows from continuous inclusion
$W\sp{2,2}(\Omega)\hookrightarrow
C\sp 0(\Omega)\hookrightarrow L\sp{2p}(\Omega),$ that
$\|w\|\sb{2p}\leq\gamma
\|w\|\sb{2,2}$, for all $w\in W\sp{2,2}(\Omega)$. From this and by
(\ref{eqn:6}) we
get (\ref{eqn:7}). Similarly, if $N=4$, (\ref{eqn:7}) follows from
$W\sp{2,2}(\Omega)\hookrightarrow L\sp{2p}(\Omega)$ together with
(\ref{eqn:6}).
Finally, in case $N>4$ and $1/2<p\leq N/(N-4)$, (\ref{eqn:7}) follows, by
using the
same argument as above, from $W\sp{2,2}(\Omega)\hookrightarrow
L\sp{\frac{2N}{N-4}}(\Omega)\hookrightarrow L\sp{2p}(\Omega)$.
\end{proof}

The assumptions on $f$ we shall assume is as follows.

\begin{itemize}
\item[$(\rm{f\sb 1})$]
$f\colon\Omega\times\mathbb{R}\times\mathbb{R}\sp{+}\to\mathbb{R}$ is a
Caratheodory
function;

\item[$(\rm{f\sb 2})$] there exists $\mu>0$ such that $N\sb{f\sb\mu}\circ
(-\Delta)\sp{-1}$ maps a ball $B_R$ of $L\sp 2(\Omega)$ into itself.
\end{itemize}

Following the same arguments as in \cite{CSBarroso} together with Theorem
\ref{trm:3.2}, we can prove the following result.

\begin{theorem}\label{trm:existence} Assume $(\rm{f\sb 1})$ and $(\rm{f\sb
2})$. Then problem (\ref{eqn:5}) has a strong solution for every
$\lambda\geq 0$.
\end{theorem}

We illustrate Theorem \ref{trm:existence} by means of the following simple
examples.

\begin{example}Consider the problem
\begin{equation}\label{eqn:elliptceq}
-\Delta u+\lambda u=\mu |u|\sp{p-2}u+a|u|\sp{q-2}u+h(x)\mbox{ in
}\Omega,\quad
u=0\mbox{ on } \partial\Omega,
\end{equation}
\noindent where $p>2$, $a\geq 0$, $3/2\leq q<2$ and $h\in L\sp 2(\Omega)$.
We now
claim that if $(H\sb{p-1})$ is fulfilled, then there exists $\lambda\sp *>0$
such
that for every $\lambda\geq -\lambda\sp *$, the problem
(\ref{eqn:elliptceq}) has at
least one strong solution. To this end, let us suppose first that
$\lambda\geq 0$.
Thanks to Theorem \ref{trm:existence}, it is enough to show that the
function
$f:\Omega\times\mathbb{R}\times\mathbb{R}\sp{+}\to\mathbb{R}$ defined by
$$
f(x,u, \mu)=\mu |u|\sp{p-2}u+a|u|\sp{q-2}u+h(x),
$$
satisfies condition $\rm{(f\sb 2)}$. From the fact that $(H\sb{p-1})$ holds
and by
(\ref{eqn:7}), we obtain the following estimate.
\begin{equation}\label{Estimate}
\|N\sb{f\sb\mu}\circ L\sp{-1}(u)\|\sb 2\leq \mu
\gamma\sp{p-1}\|u\|\sb{2}\sp{p-1}+a\gamma\sp{q-1}\|u\|\sb{2}\sp{q-1}+\|h\|\sb
2,
\end{equation}
for all $u\in L\sp 2(\Omega)$. Consequently, if $\|u\|\sb 2\leq R$ then we
get
$$
\|N\sb{f\sb\mu}\circ L\sp{-1}(u)\|\sb 2\leq
\mu\gamma\sp{p-1}R\sp{p-1}+a\gamma\sp{q-1}R\sp{q-1}+\|h\|\sb 2.
$$
Now, since $0<q-1<1$, we can choose $R>0$ large enough such that
$$
a\gamma\sp{q-1}R\sp{q-1}+\|h\|\sb 2<R.
$$
Then, choosing
$$
\mu\sp *=\frac{R-a\gamma\sp{q-1}R\sp{q-1}-\|h\|\sb
2}{\gamma\sp{p-1}R\sp{p-1}},
$$
we conclude that $\|N\sb{f\sb\mu}\circ L\sp{-1}(u)\|\sb 2\leq R$,
for all $\|u\|\sb 2\leq R$ and $\mu\in (0,\mu\sp *)$. This implies
$\rm{(f\sb 2)}$ and proves the claim for $\lambda\geq 0$. Let now
$\lambda\sp *=1/\|L\sp{-1}\|\sb{\mathcal{L}(L\sp 2(\Omega))}$ and
let $\lambda\in (-\lambda\sp *,0)$. Then $|\lambda|\cdot
\lambda\sp *<1$ and from estimate (\ref{Estimate}) we can apply
Theorem \ref{trm:3.1} to get a fixed point to (\ref{eqn:fpprb}).
Such a function will be a strong solution to
(\ref{eqn:elliptceq}). The claim is proved.
\end{example}
\par
In next example we explore the case where $\mu=1$ in  (\ref{eqn:elliptceq}).

\begin{example} Let $p$ and $h$ be as above. Consider the problem
\begin{equation}\label{pb2}
-\Delta u=|u|\sp{p-2}u+h(x)\mbox{ in  }\Omega,\qquad u=0\mbox{ on
}\partial\Omega.
\end{equation}
In this case we can obtain a strong solution to (\ref{pb2}) via
fixed point Theorem since $h$ is small enough. Indeed, define the
operator $A$ by
$$
A(u)=|L\sp{-1}(u)|\sp{p-2}L\sp{-1}(u).
$$
From inequality (\ref{eqn:7}), it follows that $A$ is well-defined
from $L\sp 2(\Omega)$ into itself. In addition, the same argument
used in (\ref{Estimate}) shows that
$$
\|A(u)\|\sb 2\leq \gamma\sp{p-1}\|u\|\sb 2\sp{p-1},
$$
for all $u\in L\sp 2(\Omega)$. In view of Theorem \ref{trm:3.4},
problem (\ref{pb2}) is solvable for $h$ small enough.
\end{example}

\begin{remark} It is worthwhile to point out that in both above
examples, the power nonlinearity $p$ might hit the critical
exponent, i.e, $p=2\sp *=\frac{2N}{N-2}$, for all $N\geq 3$.
\end{remark}
\section{A nonlinear integral equation: existence theory}
In this section we deal with the following integral equation
\begin{equation}\label{eqn:4.1}
u(t)=f(u)+\int_0\sp t g(s,u)ds,\quad u\in C(I,E),
\end{equation}
where $E$ is a reflexive space and $I=[0,T]$. Assume that the functions
involved in
equation (\ref{eqn:4.1}) satisfy the following conditions
\begin{itemize}
\item[$(H_1)$] $f:E\to E$ is sequentially weakly continuous;\\
\item[$(H_2)$] $\|f(u)-f(v)\|\leq\lambda \|u-v\|$, $(\lambda<1)$ for all
$u,v\in E$;\\
\item[$(H_3)$] $\|u\|\leq \|u-\big(f(u)-f(0)\big)\|$, for all $u\in E$;\\
\item[$(H_4)$] for each $t\in I$, the map $g\sb t=g(t,\cdot):E\to E$ is
sequentially weakly continuous;\\
\item[$(H_5)$] for each $u\in C(I,E)$, $t\mapsto g(t,u)$ is Pettis
integrable on $I$;\\
\item[$(H_6)$] there exists $\alpha\in L\sp 1[0,T]$ and  a nondecreasing
continuous function
$\varphi:[0,\infty)\to(0,\infty)$ such that $\|g(s,u)\|\leq
\alpha(s)\varphi(\|u\|)$ for a.e $s\in [0,t]$, and
all $u\in E$. Moreover, $\int\sb0\sp T\alpha(s)ds<\int\sb
0\sp\infty\frac{dx}{\varphi(x)}$.
\end{itemize}
Our existence result for (\ref{eqn:4.1}) is as follows.
\begin{theorem}\label{trm:3} Under assumptions $(H_1)$-$(H_6)$, equation
(\ref{eqn:4.1}) has at least one solution $u\in C(I,E)$.
\end{theorem}
\begin{proof} Let us define the functions
$$
J(z)=\int_{|f(0)|}\sp z \frac{dx}{\varphi(x)}\quad\mbox{ and }\quad
b(t)=J\sp{-1}\Big(\int_0\sp
t\alpha(s)ds\Big).
$$
We now define the set
\begin{eqnarray*}
M=\{u\in C(I,E):\, \quad \|u(t)\|\leq b(t)\mbox{ for all } t\in I\}.
\end{eqnarray*}

Our strategy is to apply Theorem \ref{trm:1} in order to find a fixed point
for the
operator $A+B$ in $M$, where $A,B:M\to C(I,E)$ are defined by
\begin{eqnarray*}
Au(t) &=& f(0)+ \int_0\sp tg(s,u)ds,\qquad \mbox{ and } \\
Bu(t) &=& f(u(t))-f(0).
\end{eqnarray*}
The proof will be given in several steps.

\vspace{.2cm}
\paragraph{\bf Step 1.} $M$ is bounded, closed and convex in $C(I,E)$.
\par The fact that $M$ is bounded and closed comes directly from its
definition. Let us show $M$ is convex. Let $u,v$ be
any two points in $M$. Then, there holds
\begin{eqnarray*}
&&\,\|(1-s)u(t)+sv(t)\|\leq b(t)
\end{eqnarray*}
for all $t\in I$, which implies that $(1-s)u+sv\in M$, for all $s\in[0,1]$.
This shows that $M$ is convex.

\vspace{.2cm}
\paragraph{\bf Step 2.} $A(M)\subseteq M$, $A(M)$ is weakly equicontinuous
and   $A(M)$ is relatively  weakly compact.
\par
\rm{\bf i.} Let $u\in M$ be an arbitrary point. We shall prove  $Au\in M$.
Fix $t\in I$ and consider $Au(t)$.
Without loss of generality, we may assume that $Au(t)\neq 0$. By the
Hahn-Banach Theorem there exists $\psi_t\in
E\sp*$ with $\|\psi_t\|=1$ such that $\langle\psi_t,
Au(t)\rangle=\|Au(t)\|$. Thus,
\begin{eqnarray}\label{eqn:4.2}
\|Au(t)\|&=&\langle\psi_t, Au(t)\rangle=\langle\psi_t,f(0)\rangle+\int_0\sp
t\langle \psi_t, g(s,u)\rangle ds\\
&\leq& \|f(0)\|+\int_0\sp t\alpha(s)\varphi(\|u(s)\|)ds\nonumber \\
&\leq&\|f(0)\|+\int_0\sp t\alpha(s)\varphi(b(s))ds\leq b(t),\nonumber
\end{eqnarray}
since
$$
\int_{|f(0)|}\sp{b(s)}\frac{dx}{\varphi(x)}=\int_0\sp s\alpha(x)dx.
$$
Then, $(\ref{eqn:4.2})$ implies that $A(M)\subseteq M$. Analogously one
shows that,
\begin{eqnarray}\label{eqn:4.3}
\|Au(t)-Au(s)\|&\leq& \int_t\sp s
\alpha(\eta)\varphi(\|u(\eta)\|)d\eta\nonumber\\
&\leq& \int_t\sp s\alpha(\eta)\varphi(b(\eta))d\eta=\int_t\sp s
b'(\eta)d\eta\nonumber\\
&\leq &|b(t)-b(s)|,
\end{eqnarray}
for all $t,s\in I$. Thus it follows from (\ref{eqn:4.2}) that $A(M)$ is
weakly
equicontinuous.
\par
\rm{\bf ii.} Let $(Au\sb n)$ be any sequence in $A(M)$. Notice
that $M$ is bounded. By reflexiveness, for each $t\in I$ the set
$\{Au\sb n(t):\, n\in\mathbb{N}\}$ is relatively weakly compact.
As before, one shows that $\{Au\sb n:\, n\in\mathbb{N}\}$ is a
weakly equicontinuous subset of $C(I,E)$. It follows now from the
Ascoli-Arzela Theorem that $(Au\sb n)$ is relatively weakly
compact, which proves the third assertion of $\rm{Step~ 2}$.

\vspace{.2cm}
\paragraph{\bf Step 3.} $A$ is sequentially weakly continuous.
\par
Let $(u\sb n)$ be a sequence in $M$ such that $u\sb
n\rightharpoonup u$ in $C(I,E)$, for some $u\in M$. Then, $u\sb
n(s)\rightharpoonup u(s)$ in $E$ for all $s\in I$. By assumption
$(H_5)$ one has that $g(s,u\sb n(s))\rightharpoonup g(s,u(s))$ in
$E$ for all $s\in I$. The Lebesgue dominated convergence Theorem
yields that $Au\sb n(t)\rightharpoonup Au(t)$ in $E$ for all $t\in
I$. On the other hand, it follows from (\ref{eqn:4.3}) that the
set $\{Au\sb n:\, n\in\mathbb{N}\}$ is a weakly equicontinuous
subset of $C(I,E)$. Hence, by the Ascoli-Arzela Theorem there
exists a subsequence $(u\sb{n_j})$ of $(u_n)$ such that
$Au\sb{n_j}\rightharpoonup v$ for some $v\in C(I,E)$.
Consequently, we have that $v(t)=Au(t)$ for all $t\in I$ and hence
$Au\sb{n_j}\rightharpoonup Au$. Now, a standard argument shows
that $Au\sb n\rightharpoonup Au$. This proves $\rm{Step~ 3}$.

\vspace{.2cm}
\paragraph{\bf Step 4.} $B$ satisfies conditions $(b)$ and $(d)$ of
Theorem
\ref{trm:1}.\par By $(H_2)$ clearly we see that $B$ is a
$\lambda$-contraction in
$C(I,E)$. Now, in order to verify condition $(d)$ to $B$, we first remark
that by
combining (\ref{eqn:4.3}) with $(H_2)$, it follows that
  $\mathcal{F}$ is weakly equicontinuous in $C(I,E)$. So is
$B(\mathcal{F})$. Let now
  $(u\sb n)\subset\mathcal{F}$ be such that $u\sb n\rightharpoonup u$, for
some $u\in M$. Then by assumption
   $(H_1)$, we obtain $Bu\sb n (t)\rightharpoonup Bu(t)$. Since $(Bu\sb n)$
is weakly equicontinuous in $C(I,E)$ and
   $\|(Bu\sb n)(t)\|\leq \lambda \|u\sb n(t)\|$ holds for all
$n\in\mathbb{N}$, we may apply the Ascoli-Arzela Theorem and
concludes that there exists a subsequence $(u\sb{n_j})$ of $(u_n)$
such that $Bu\sb{n_j}\rightharpoonup v$, for some $v\in C(I,E)$.
Hence, $Bu=v$ and by standard arguments we have $Bu\sb
n\rightharpoonup Bu$ in $C(I,E)$. This completes Step 4.

\vspace{.2cm}
\paragraph{\bf Step 5.} Condition $(c)$ of Theorem \ref{trm:1} holds.
\par
Suppose that $u=Bu+Av$ for some $v\in M$. We will show that $u\in M$. By
condition $(H_3)$ it follows that
\begin{eqnarray*}
\|u(t)\|\leq \|u(t)- Bu(t)\|=\|Av(t)\|.
\end{eqnarray*}
Once $v\in M$ implies $Av\in M$, we conclude $u\in M$.

\vspace{.2cm}
\paragraph{\bf Step 6.} Condition $(e)$ of Theorem \ref{trm:1} holds.
\par
Let $(u\sb n)\subset\mathcal{F}$ be an arbitrary sequence. Then,
$(u\sb n)$ is weakly equicontinuous in $C(I,E)$. Also, one has
that
$$
\|u\sb n(t)\|\leq (1-\lambda)\sp{-1}\cdot b(t),
$$
for all $t\in I$, that is, for each $t\in I$ the set $\{u\sb
n(t)\}$ is relatively weakly compact in $E$. Thus, invoking again
the Ascoli-Arzela Theorem we obtain a subsequence of $(u\sb n)$
which converges weakly in $C(I,E)$. By the Eberlein-\v{S}mulian
Theorem, it follows that $\mathcal{F}$ is relatively weakly
compact.
\par
Theorem \ref{trm:1} now gives a fixed point for $A+B$ in $M$, and
hence a solution to (\ref{eqn:4.1}). Such a solution is, a priori,
a weakly continuous curve; however since $u \in
\overline{co}(\mathcal{F})$ we can conclude $u$ is actually norm
continuous.
\end{proof}
\begin{remark} Theorem \ref{trm:3} is a generalization of Theorem 2.4 in
\cite{6}.
\end{remark}

We complete this section by presenting a wide and illustrative
class of maps $f$ defined on the real line fulfilling condition
$(H\sb 3)$ in Theorem \ref{trm:3}. Let $ f \colon \mathbb{R} \to
\mathbb{R}$ satisfy
$$
    tf(t) \le 0.
$$
It is easy to check that functions satisfying the above inequality
fulfill assumption $(H\sb 3)$ in Theorem \ref{trm:3}.


\section{Optimal regularity for equation $(\ref{eqn:4.1})$ via topological methods}

In this section we want to explore the optimal regularity of equation
$(3.1)$.
\begin{equation} \label{IntEq}
    u(t) = f(u(t)) + \int_0^t g(s,u(s))ds
\end{equation}

Questions related to regularity of solutions to abstract nonlinear
integral equations is quite important. The importance of this
question is intrinsically related to almost everywhere
differentiability of curves in Banach spaces. By knowing the
solution of a nonlinear integral equation is sufficiently regular
we can recover the original differential equation the nonlinear
integral equation models. \\
\indent The conditions on equation (\ref{IntEq}) will be slight
different. However these will not be more restrictive in practical
applications. We shall assume
\begin{itemize}
\item[$(C_1)$] $f:E\to E$ is sequentially weakly continuous;\\
\item[$(C_2)$] $f$ is 1-Lipschitz and differentiable;\\
\item[$(C_3)$] $\|u\|\leq \|u-\big(f(u)-f(0)\big)\|$, for all $u\in E$
$\&~ \|w\| \le \| w - Df(u)w \|, ~\forall w,u \in E$;\\
\item[$(C_4)$] $g$ is weak Caratheodory; \\
\item[$(C_5)$] $g$ has polynomial growth, i.e., $\|g(s,u)\| \le C\cdot(
\|u\|^r + 1)$, for some $r \ge 0$; \\
\item[$(C_6)$] There exists $R > \|f(0)\|$ such that $ C \cdot ( T^r + 1) =
\dfrac{R - \|f(0)\|}{R^r + 1}$.
\end{itemize}
\begin{remark} The differentiability asked in condition $(C_2)$ is
related to one of the most important problems in nonlinear
functional analysis: the ''almost everywhere" differentiability of
Lipschitz maps between Banach spaces. Many progress have been done
in the direction of defining the right notion of ''almost
everywhere" for Banach spaces. (see for instance \cite{JLPS} and
\cite{ Lindenstrauss_Preiss}). We should point out though that
once the space we are working on is separable and reflexive it
follows from a result due to Aronszajn, Christensen and Mankiewicz
that Lipschitz functions are Gateaux differentiable outside an
Aronszajn null set. See, for instance, Theorem $6.42$ in \cite{BL}. \\
\indent Condition $(C_6)$ is always verified by taking $T$ small
enough. Such a condition is actually a restriction on the global
definition of solutions. Such a constraint appears even in very
simple versions of equation (\ref{IntEq}). It is possible though,
to show that if $r\le 1$ we can find a global solution to equation
(\ref{IntEq}).
\end{remark}
\begin{theorem}[Optimal Regularity] \label{optimal} Under
assumptions $(C_1)-(C_6)$ equation (\ref{IntEq}) has at least a Lipschitz
solution.
Such a solution is almost everywhere differentiable.
\end{theorem}
Before proving Theorem \ref{optimal} we need to develop some
tools. We will start by defining a new locally convex topology to
$W^{1,\infty}(I, E)$.
\begin{definition} Let $E$ be a Banach space and $I$ a bounded interval in
$\Bbb{R}$. Let us denote by $\mathcal{T}^n$ the family of seminoms
given by $\mathcal{T}^n := \{ \rho\colon W^{1,\infty}(I,E) \to
\Bbb{R}_+ : \rho = |f| \text{ and } f \in [W^{1,n}(I,E)]^{*} \}$.
We define the T-topology in $W^{1,\infty}(I,E)$ to be the locally
convex topology generated by $\bigcup\limits_{n\ge 2}
\mathcal{T}^n$.
\end{definition}
\begin{proposition} Let $E$ be a separable reflexive Banach space. Then
$$\big( B_{W^{1,\infty}(I,E)}(R), \textrm{T} \big )$$
is metrizable.
\end{proposition}
\begin{proof}
Initially we note that for each $n\ge 2$, the space
$[W^{1,n}(I,E)]^{*}$ is separable. The Lemma follows easily from
this general metrization Theorem: {\it A topological space is
metrizable if and only if it is regular and has a basis that is
the union of at most countably many locally finite systems of open
sets. }
\end{proof}
\begin{proposition} \label{comp} Let $E$ be a separable  reflexive Banach
space. Then $$X :=\big( \overline{B_{W^{1,\infty}(I,E)}(R)}, \textrm{T} \big
)$$ is a
compact space.
\end{proposition}
\begin{proof} Let $\{u_j\} \subset X$. By the continuous embedding
$W^{1,\infty}(I,E)\hookrightarrow W^{1,n}(I,E)$ we have that $\|
u_j\|_{W^{1,n}} \le
|I|^{1/n} R$. Since $\big ( \overline{B_{W^{1,n}(I,E)}(|I|^{1/n} R)},
\mathcal{T}^n
\big )$ is a compact space, there exists a subsequence $\{u_{j_k}\}$ which
converges
weakly to a $u \in W^{1,n}(I,E)$. By the Cantor Diagonal Argument, we build
a
subsequence $\{u_{j_r}\}$ which converges weakly to $u$ in $W^{1,n}(I,E)$
for all
$n\ge 2$. This initially implies that $u_{j_r} \stackrel{\textrm{T}}{\to}
u$. In
addition $ \| u\|_{W^{1,n}} \le |I|^{1/n} R. $ Letting $n\to \infty$ leads
$u\in
W^{1,\infty}$ and $\|u\|_{W^{1,\infty}} \le R$.
\end{proof}
\begin{lemma} \label{MT} Let $f: E \to E$ be a differentiable
1-Lipschitz map. Then for each $p> 1$ and $\lambda \in (0,1)$ the map $ u
\mapsto I -
\lambda f(u)$ is a sequentially weak continuous homeomorphism between
$W^{1,p}(I,E)$
onto itself.
\end{lemma}
\begin{proof} Let $\psi \in W^{1,p}(I,E)$ be given. By changing $f$ by $ f -
f(0)$,
we may assume, without loss of generality that $f(0) = 0$, Let us start by
estimating
$\| f(\xi)\|_{W^{1,p}(I,E)}$, for any $\xi \in W^{1,p}(I,E)$:
\begin{equation} \label{estimateW}
\begin{array}{lll}

\| f(\xi)\|_{W^{1,p}} &= & \|f(\xi (t))\|_{L^p} + \|\partial_tf(\xi
(t))\|_{L^p}\\
                  &\le& \|\xi\|_{L^p}  + \|(D
f)(\xi(t))\cdot\partial_t\xi(t)\|_{L^p} \\
                  & \le & \|\xi\|_{W^{1,p}}.
\end{array}
\end{equation}
Inequality (\ref{estimateW}) tells us the substitution operator
$N_f:W^{1,p}(I, E) \to W^{1,p}(I,E)$ given by $N_f(\xi) := f(\xi)$
is a bounded operator. We claim $N_f$ is a sequentially weakly
continuous map. Indeed, Let $\{u_n\} \subset W^{1,p}(I,E)$ be such
that $u_n \rightharpoonup u$ in $W^{1,p}(I,E)$. We know that
$\|u_n\|_{W^{1,p}}$ is bounded. Hence the sequence $\{N_f(u_n)\}$
is bounded in $W^{1,p}(I,E)$. Due to the reflexibility of
$W^{1,p}(I,E)$, we may assume $N_f(u_n) \rightharpoonup v$ in
$W^{1,p}(I,E)$, for some $v\in W^{1,p}(I,E)$. The idea now is to
show that $v = A(u)$. The crucial information here is the
continuous embedding
\begin{equation} \label{embed0}
W^{1,p}(I,E)\hookrightarrow C^\alpha(I,E), \textrm{ for each } 0\le \alpha <
1-1/p.
\end{equation}
Thus we can conclude that for each $s\in I$, $u_n(s) \rightharpoonup u(s)$
in $E$.
Indeed, for a fixed $\Psi \in E^{*}$, the map $\Psi_s\colon C(I,E) \to
\Bbb{R}$
defined by $\Psi_{s} (u) := \Psi(u(s))$ is a continuous linear functional.
Therefore
we also have that
$$
    N_f(u_n)(s) \rightharpoonup N_f(u)(s), \forall s \in I.
$$
However, the same argument as before works to show that
$$
N_f(u_n)(t) \rightharpoonup v(t) \textrm{ in } E,~\forall t \in I.
$$
Thus, $v(t)$ has to be equal to $N_f(u)(t)$ for every $t\in I$. \\
Let us define $ \Lambda : W^{1,p}(\Omega) \to W^{1,p}(\Omega)$ by
$$
\Lambda(\xi) = \psi + \lambda f(\xi).
$$
We observe that once $N_f$ is sequentially weakly continuous, so is
$\Lambda$.
Moreover, the solvability of the equation
\begin{equation} \label{onto}
    (I - \lambda N_f)(u) = \psi
\end{equation}
is equivalent to finding a fixed point for $\Lambda$. For each $\xi \in
W^{1,p}(\Omega)$, we have from the triangular inequality and inequality
(\ref{estimateW})
$$
\begin{array}{lll}
  \| \Lambda(\xi) \|_{W^{1,p}} &\le& \| \psi\|_{W^{1,p}} + \lambda
  \| N_f(\xi) \|_{W^{1,p}} \\
  &\le&  \| \psi\|_{W^{1,p}} + \lambda \|\xi\|_{W^{1,p}}.
\end{array}
$$
Let us fix $ M > \dfrac{ \|\psi\|_{W^{1,p}} }{1- \lambda}$. For such an $M$
we see
that if $ \|\xi\|_{W^{1,p}(\Omega) } \le M $, then
\begin{equation} \label{LambdaBound}
     \|\Lambda(\xi)\|_{W^{1,p}} \le \| \psi\|_{W^{1,p}} + \lambda M \le M.
\end{equation}
In other words, $\Lambda$ maps the closed ball of radius $M$ in
$W^{1,p}(I,E)$ into itself. Let $X$ denote
$\overline{B_{W^{1,p}}(M)}$ endowed with the weak topology. So $X$
is a compact and convex  set of a locally convex space. In
additional, as we pointed out before, ${\Lambda} : X \to X$ is a
continuous map. Invoking the Leray-Schauder-Tychonoff fixed point
Theorem we conclude that ${\Lambda}$ has a fixed point which is
precisely a solution to (\ref{onto}). Since $\psi \in
W^{1,p}(I,E)$ was
taken arbitrarily we have proven $(I - \lambda N_f)$ is onto. \\
We now turn our attention to uniqueness. Let us suppose that there exist
$u_1,~u_2
\in W^{1,p}(I,E)$ such that
$$ (I - \lambda N_f)u_i(t) = \psi(t) \textrm{ for } i=1,2$$
Subtracting these above equations, we find
$$
f(u_1(t)) - f(u_2(t)) = \lambda (u_1(t) - u_2(t)).
$$
Therefore
$$ \|f(u_1(t)) - f(u_2(t))\| = \lambda \|u_1(t) - u_2(t)\| \le \|u_1(t) -
u_2(t)\|. $$
If $ \|u_1 - u_2\| > 0$ in a set of positive measure, we would
find, $\lambda \ge 1$. Hence the solution of $(P)$ is unique. \\
\indent Finally let us study the weak sequential continuity of $R_\lambda =
(I -
\lambda N_f)^{-1}:W^{1,p}(I,E) \to W^{1,p}(I,E)$. Suppose $R_\lambda (\psi)
= u$,
i.e., $ (I - \lambda N_f)u = \psi $. Then
$$
\begin{array}{lll}
\|\psi\|_{W^{1,p}} &\ge& \|u\|_{W^{1,p}} - \lambda \|f(u)\|_{W^{1,p}} \\
                   &\ge& (1-\lambda)\|u\|_{W^{1,p}}.
\end{array}
$$
Writing in a better way, $\|R_\lambda(\psi)\|_{W^{1,p}} \le \dfrac{
\|\psi\|_{W^{1,p}} }{1-  \lambda } $.  We have just verified that
$R_\lambda$ is a
bounded operator. Suppose $\psi_n \rightharpoonup \psi$ in $W^{1,p}(I,E)$.
Let us
denote by $u_n = R_\lambda(\psi_n)$. The sequence $\{u_n\}$ is bounded,
therefore, up
to a subsequence, we may assume that $ u_n \rightharpoonup u$ in
$W^{1,p}(I,E)$. From
the weak sequential continuity of $(I - \lambda N_f)$ we have
$$
\psi_n =  (I - \lambda N_f)(u_n) \rightharpoonup (I - \lambda N_f)(u).
$$
This implies that $R_\lambda(\psi) = u$, and thus, $R_\lambda(\psi_n)
\rightharpoonup
R_\lambda(\psi)$ as desired.
\end{proof}
\noindent {\it Proof of Theorem \ref{optimal}. } Define $ M := \{ u \in
W^{1,\infty}(I,E) : \| u\|_{1,\infty} \le R \}$, where $R$ is given by
condition
$(C_6)$. Define also,
$$
A(u) := f(0) + \int_0^t g(s,u(s))ds
$$
$$
B(u) = f(u) - f(0)
$$
The strategy is to find a fixed point in $M$ to the operator
$A+B$. \\
{\bf Step 1}. $A$ maps $W^{1,p}(I,E)$ into $W^{1, p/r}(I,E)$.
Moreover it is sequentially weakly continuous.\\
Indeed, Let us estimate $\|A(u)\|_{W^{1,p/r}}$. We first deal with
$\|A(u)\|_{L^{p/r}}$.
$$
\|A(u)\|_{L^{p/r}} \le \|f(0)\|_E \cdot |I|^{r/p} +\left \{ \int_0^T \left
\|
\int_0^t g(s,u(s))ds \right \|^{p/r}_E dt \right \}^{r/p}.
$$
By Jensen's inequality we obtain
$$
\begin{array}{lll}
\left \{ \displaystyle\int_0^T \left \| \int_0^t g(s,u(s))ds \right
\|^{p/r}_E dt
\right \}^{r/p} &\le& \left \{ \displaystyle\int_0^T
t^{p-1} \int_0^t \|g(s,u(s) \|^{p/r}_E ds ~dt \right \}^{r/p} \\
& \le & \left \{ T^p \displaystyle\int_0^T C^{p/r}\big( \|u(s)\|^{r}_E + 1
\big )^{p/r} ds \right \}^{r/p} \\
& \le & T^r C \big (  \|u\|_{L^p}^r +  |I|^{r/p} \big ).
\end{array}
$$
We have shown that
\begin{equation} \label{Lp}
\| A(u) \|_{L^{p/r}} \le \|f(0)\|_E \cdot T^{r/p} + {T^r} C\big (
\|u\|_{L^p}^r +
T^{r/p} \big ).
\end{equation}
In the same way we find
\begin{equation} \label{D-Lp}
\| \partial_t A(u) \|_{L^{p/r}} \le {C ( \|u\|_{L^p}^r + T^{r/p}) }.
\end{equation}
Adding up the above two inequalities we end up with the following inequality
\begin{equation} \label{W1,p}
\| A(u) \|_{W^{1,p/r}} \le {C (T^r +1)}\|u\|_{L^p}^r + T^{r/p} \left (
\|f(0)\|_E  +
({T^r+1})C  \right ).
\end{equation}
Let us now show he map $A\colon  W^{1,p}(I,E) \to W^{1,p/r}(I,E)$ is
sequentially
weakly continuous. Let $\{u_n\} \subset W^{1,p}(I,E)$ be such that $u_n
\rightharpoonup u$ in $W^{1,p}(I,E)$. Immediately we know that
$\|u_n\|_{W^{1,p}}$
and $\sup_{n,s} \|u_n(s)\|_E$ are bounded. The last bound is due to the
continuous
embedding,
\begin{equation} \label{embed}
W^{1,p}(I,E)\hookrightarrow C^\alpha(I,E), \textrm{ for each } 0\le \alpha <
1-1/p.
\end{equation}
Inequality (\ref{W1,p}) says in particular that $A$ is a bounded
operator. Hence the sequence $\{A(u_n)\}$ is bounded in
$W^{1,p/r}(I,E)$. Due to the reflexibility of $W^{1,p/r}(I,E)$, we
may assume $A(u_n) \rightharpoonup v$ in $W^{1,p/r}(I,E)$, for
some $v\in W^{1,p/r}(I,E)$. The idea now is to show that $v =
A(u)$.  We remark, as done in the proof of Lemma \ref{MT}, that
for each $s\in I$, $u_n(s) \rightharpoonup u(s)$ in $E$.  Let us
fix a $\Phi \in E^{*}$. For each $s\in I$ which the map
$g(s,\cdot)\colon E \to E$ is sequentially weakly continuous,
there holds
$$
\Phi(g(s,u_n(s))) \longrightarrow \Phi(g(s,u(s))).
$$
Moreover,
$$
\begin{array}{lll}
|\Phi(g(s,u_n(s)))| & \le & \| \Phi \|_{E^{*}} \cdot \| g(s,u_n(s)) \|_E
\\
& \le & \| \Phi \|_{E^{*}} \big ( C\|u_n(s)\|_E^r + A \big ) \\
& \le & \| \Phi \|_{E^{*}}\cdot \widetilde{C}.
\end{array}
$$
Follows now from Lebesgue's dominated convergence Theorem that
$$
\langle ~\Phi, A(u_n)(t) ~\rangle := \Phi(f(0)) + \int_0^t
\Phi(g(s,u_n(s)))ds
\stackrel{n\to\infty}{\longrightarrow} \langle ~\Phi, A(u)(t) ~\rangle.
$$
We have proven that for every $t\in I$, $ A(u_n)(t) \rightharpoonup A(u)(t)
$ in $E$.
However, as we did before, the continuous embedding,
$W^{1,p/r}(I,E)\hookrightarrow
C^\gamma(I,E)$, for all $0\le \gamma < 1- r/p$, yields
$$
A(u_n)(t) \rightharpoonup v(t) \textrm{ in } E.
$$
Thus, $v(t)$ has to be equal to $A(u)(t)$ for every $t\in I$. \\
{\bf Step 2}. $A$ maps $M$ into itself.\\
By letting $p \to \infty$ in (\ref{W1,p}) we obtain that $A\colon
W^{1,\infty}(I,E)
\to W^{1,\infty}(I,E)$ and
\begin{equation} \label{W1,infty}
\| A(u) \|_{W^{1,\infty}} \le {C (T^r+1)}\cdot (\|u\|^r_{L^{\infty}} + 1) +
\|f(0)\|_E.
\end{equation}
Notice that if $ \| u \|_{W^{1,\infty}} \le R$ there holds
$$
\begin{array}{lll}
\| A(u) \|_{W^{1,\infty}} &\le& {C (T^r+1)}\cdot (R^r+
1) + \|f(0)\|_E \\
& = & R.
\end{array}
$$
{\bf Step 3}. The map $B$ is sequentially weakly continuous from
$W^{1,p}(I,E)$ into itself. \\
Indeed we estimate
$$
\| B(u)\|_p^p = \int_0^T \| f(u(t)) - f(0)\|^pdt \le \int_0^T \|u(t)\|^p.
$$
$$
\| \partial_t B(u)\|_p^p = \int_0^T \| Df(u(t))u_t(t)\|^pdt \le \int_0^T
\|u_t(t)\|^p.
$$
The above inequalities show $B$ is a bounded operator from $W^{1,p}(I,E)$
into
itself. We now argue as in Step 1 to conclude
the sequential weak continuity of $B$. \\
{\bf Step 4}. $M$ is $T$-compact and $A,B \colon M \to
W^{1,\infty}(I,E)$ are $T$-continuous.\\
The fact that $M$ is $T$-compact follows from proposition \ref{comp}. Let us
show $A$
is $T$-continuous: Let $\{u_j\} \subset M$ be such that $ u_j
\stackrel{\textrm{T}}{\to} u$. It means $u_j \rightharpoonup u$ in
$W^{1,m}(I,E)$,
for each $ m\ge 2$. We have to show that $ A(u_j) \rightharpoonup A(u) $ in
$W^{1,n}(I,E)$ for every $n\ge 2$. To this end, let us fix $n\ge 2$ and let
$m =
\lceil n\cdot r \rceil$, the lowest integer bigger than $ n \cdot r$. We
know $u_j
\rightharpoonup u$ in $W^{1,m}(I,E)$. By Step 1, $A(u_j) \rightharpoonup
A(u)$ in
$W^{1,m/r}(I,E)$. Once $W^{1,m/r}(I,E)$ is continuously embedded into
$W^{1,n}(I,E)$,
the convergence also holds in $W^{1,n}(I,E)$. We have proven $A$ is
$T$-continuous. A
similar argument shows $B$
is also $T$-continuous.\\
{\bf Step 5}. For all $\lambda \in (0,1)$, $(I- \lambda B)$ is
$T$-homeomorphism from $W^{1,\infty}(I,E)$ onto itself. \\
Let $\psi \in W^{1,\infty}(I,E)$ be given. For each $p\ge 1$, we
know $\psi \in W^{1,p}(I,E)$. We then can apply Lemma \ref{MT} to
conclude there is a unique $u \in W^{1,p}(I,E)$, which solves
$$
    (I- \lambda B)(u) = \psi.
$$
Moreover, it follows from inequality (\ref{LambdaBound}) that
$$
    \| u \|_{W^{1,p}} \le \dfrac{\|\psi\|_{W^{1,p}}}{1- \lambda} +
    \varepsilon, ~\forall  \varepsilon > 0.
$$
Letting $\varepsilon \to 0$ and $p \to \infty$ we conclude $u$ in
in $W^{1,\infty}(I,E)$ and moreover $\|u\|_{W^{1,\infty}} \le
\frac{\|\psi\|_{W^{1,\infty}}}{1-\lambda}$. Lemma \ref{MT} also
provides the sequential weak continuity of $(I - \lambda B), ~(I -
\lambda B)^{-1} \colon W^{1,p}(I,E) \to W^{1,p}(I,E)$. This
together with an argument like in Step 4, concludes the proof of Step 5.\\
{\bf Step 6}. Condition $(c)$ of main Corollary \ref{cor:1} holds.\\
Indeed, suppose
$$
    u = \lambda B(u) + A(v), \textrm{ for some } v \in M.
$$
From condition $(C_3)$,
$$
\| u(t) \| \le \| u(t) - [ f(u(t)) - f(0) ] \| = \| A(v(t)\| \le R.
$$
Analogously,
$$
\| u_t(t) \| \le \|\big ( I - Df(u(t)) \big )u_t(t) \| = \| g(t, v(t)) \|
\le R.
$$
Thus $\| u\|_{W^{1,\infty}} \le R$. \\
We have verified all the hypothesis of Corollary \ref{cor:1},
which assures a fixed point to $A+B$. Such a fixed point is
Lipschitz since it lies in $W^{1,\infty}(I,E)$. Furthermore, since
$E$ is reflexive it has the Radon-Nikod\'ym property. Therefore
the Lipschitz solution $u\colon I \to E$ is almost everywhere
differentiable. \hfill $\Box$

\begin{remark} When we are dealing with real valued functions
there is a simple way of verifying condition $(C_3)$ in Theorem
\ref{optimal}. Indeed if $f$ is nonincreasing and satisfies $tf(t)
\le 0$, it is easy to verify that condition $(C_3)$ is fulfilled.
\end{remark}


\section{A geometric approach to fixed point results for sum of operators on uniformly convex
spaces.} In this section we are interested in a variant of
equation (\ref{eqn:4.1}). Indeed we shall work on the following
nonlinear integral equation:
\begin{equation}
\label{eqn:4}
    u(t)=f(t,u(t))+\Phi \left (t, \int_0\sp t k(t,s)u(s)ds \right ),\quad
u\in L\sp p(I,E),
\end{equation}
where $E$ is a uniformly convex space, $1< p < \infty$ and
$I=[0,T]$.  The strategy here is rather different from the
strategy used in the previous sections. Indeed we shall explore a
geometric idea to guarantee somehow the hard-to-check assumption
of Theorem \ref{trm:I}, i.e., condition $(c)$. We shall develop
these ideas directly to analyze the important nonlinear integral
equation (\ref{eqn:4}). \\
\indent Let us point out that, a priori, equation (\ref{eqn:4}) is
more delicate than (\ref{eqn:4.1}), since the nonlinearity $\Phi$
nulls the regularity property of the integral. In the study of
equation (\ref{eqn:4}), the natural assumptions on $f,~ k$ and
$\Phi$ are:
\begin{itemize}

\item[$(A_1)$] $f\colon I \times E\to E$ is a measurable family of maps
satisfying $\| f(t,x) - f(t,0) \| \le \|x\| ~
\forall x \in E.$ \\

\item[$(A_2)$] $k \in L^\infty(I,L^q(0,T))$, where $\frac{1}{p} +
\frac{1}{q} = 1$. We shall denote by $C := \| k
\|_\infty$.\\

\item[$(A_3)$] $\Phi$ \textrm{ is a weak Carath\'eodory map satisfying }

    $\|\Phi(t,u)\|_{E} \le G(t)\psi(\|u\|_{E})$, \textrm{ where,}
\begin{equation} \label{A_4}
    \left \{
    \begin{array}{l}
    G \in L^p(0,T) \\
    \psi \in L^\infty_{\mathrm{loc}}(0,T) \textrm{ and } \exists
\overline{R} \textrm{ with } \dfrac{\|G\|_p\cdot
    \psi( C\cdot \overline{R})}{\overline{R} - \|f(t,0)\|_{p}} \le 1.
    \end{array} \right.
\end{equation}
\end{itemize}

\indent Let us recall that a map $\Phi \colon I \times E \to E$ is said to
be a weak Carath\'eodory map if for
each $x \in E$ fixed, the map $ t \mapsto \Phi(t,x)$ is measurable and for
almost every $t\in I$ the map $ x
\mapsto \Phi(t,x)$ is sequentially weakly continuous.
\begin{remark} As we shall see in the proof of Theorem \ref{Existence},
condition $(A_1)$ only need to be held
for $ x \in B_E(\overline{R})$.
In assumption ($A_3$) we may always assume $\psi$ nondecreasing and
everywhere defined. The existence
  of such a $\overline{R}$ in hypothesis (\ref{A_4}) is less restrictive
than natural assumptions on the asymptotic
  behavior of $\psi$.
\end{remark}
 To grapple with the difficulty of losing the regularity properties of the integral,
 we will need to develop a sequential weak continuity
result for a class of operators acting on vector-valued Lebesgue
spaces. Moreover a geometric assumption will also be needed to
assure the existence of a solution to problem (\ref{eqn:4}). Such
an assumption is somehow related to monotonicity hypothesis on the
operators involved in problem (\ref{eqn:4}). Let us start by the
definitions and main results involved in such a geometric
condition.
\begin{definition} \label{angle} Let $E$ be a normed vector space. We
define the notion of angle between two
  nonzero vectors $x,~y$ as follows:
  $$
    \alpha(x,y) := \left \| \dfrac{x}{\|x\|_E} - \dfrac{y}{\|y\|_E} \right
\|_{E}
  $$
\end{definition}
Let now $E$ be a uniformly convex space. Its modulus of convexity, $\delta$,
is defined as
$$
\sup \left \{ \left \| \dfrac{x+y}{2} \right \|_E : \|x\|_E = \|y\|_E = 1;~
\| x - y \|_E = \varepsilon \right
\} = 1 - \delta(\varepsilon).
$$
\begin{lemma} \label{Triang} Let $v_1,~v_2,~...,~ v_n$ be nonzero elements
of a uniformly convex space $E$.
Suppose $ V := \sum_{i=1}^n v_i \not = 0$. Let us denote by $ \alpha_i =
\alpha( v_i, V)$. Then
$$
    \| V \|_{E} \le \sum\limits_{i=1}^n \Big ( 1 - 2\delta(\alpha_i) \Big )
\cdot \|v_i\|_E,
$$
where $\delta$ is the modulus of convexity of $E$.
\end{lemma}
\begin{proof} It follows from the definition of the modulus of convexity
that for each $i$, running from $1$ to
$n$, we have
$$
    \Big \| \|V\|_E v_i + \|v_i\|_E V \Big \|_E \le 2\big ( 1 -
\delta(\alpha_i) \big ) \|V\|_E\cdot \|v_i\|_E.
$$
Summing the above inequality over $i$ we find
$$
\sum\limits_{i=1}^n \Big \| \|V\|_E v_i + \|v_i\|_E V \Big \|_E \le 2
\sum\limits_{i=1}^n\big ( 1 -
\delta(\alpha_i) \big ) \|V\|_E\cdot \|v_i\|_E.
$$
We now apply the standard triangular inequality to the left hand side of the
above inequality and end up with
$$
\| V \|_E \cdot \Big ( \| V \|_E + \sum\limits_{i=1}^n \|v_i\|_E \Big ) \le
2 \|V\|_E \sum\limits_{i=1}^n\big (
1 - \delta(\alpha_i) \big ) \|v_i\|_E.
$$
Cancelling $\|V\|_E$ out from the above and rearranging the
reminder part we conclude the Lemma.
\end{proof}
\begin{definition} Let $E$ be a uniformly convex space. We define
$\epsilon_0 = \epsilon_0(E) > 0$ to be the smallest
positive number such that whenever we write $ \epsilon_0 = \epsilon_1 +
\epsilon_2 $ with $0 \le \epsilon_1,
\epsilon_2 \le 2$, we have
\begin{equation}
\delta( \epsilon_1) + \delta (\epsilon_2) \ge 1/2.
\end{equation}
\end{definition}
\begin{definition} Let $C$ be a cone in a normed vector space. We define its
opening as follows:
$$
\theta(C) : = \sup \{ \alpha(x,y) : x,y \in C \}.
$$
\end{definition}
Let us relate these definitions to problem (\ref{eqn:4}). We define  $g(t,x)
:= f(t,x) - f(t,0)$ and $B \colon
L^p(I,E) \to L^p(I,E)$ by
$$
    B(u)(t) = g(t,u(t)).
$$
Also we define $A \colon L^p(I,E) \to L^p(I,E)$ by
$$
A (u)(t):= f(t,0) + \Phi \left (t, \int_0\sp t k(t,s)u(s)ds \right ).
$$
The condition on $B$ is:
\begin{itemize}
\item[$(A_4)$]  $B \colon L^p(I,E) \to L^p(I,E)$ is sequentially weakly
continuous.
\end{itemize}
\begin{remark} The sequential weak continuity of substitution operators
acting on vector-valued Lebesgue
spaces was studied in \cite{MT} and \cite{MT2}. We should mention that in
many practical applications which
arise from physical models one can verify condition $(A_4)$. It has been of
large interest of the authors the
study of sufficient condition to assure condition $(A_4)$. One of the
simplest technic which has contemplated
many practical situations is the following: There exists a Banach space $F$
such that $E$ is compacted embedded
into $F$ and for some $s<0$ the operator $B \colon L^p(I,E) \to L^p(I,E)$
extends to $B \colon W^{s,p}(I,F) \to
W^{s,p}(I,F)$ in a demicontinuous fashion, where $W^{s,p}(I,F)$ is the
Sobolev-Slobodeckii spaces. Indeed if the
above holds, $L^p(I,E)$ is compacted embedded into $W^{s,p}(I,F)$ (see
\cite{Amann}). Let $u_n \rightharpoonup
u$ in $L^p(I,E)$. Since $B$ is a bounded operator, up to a subsequence we
might assume that $B(u_n)
\rightharpoonup v$ for some $v$ in $L^p(I,E)$. By the compact embed
$L^p(I,E) \hookrightarrow W^{s,p}(I,F)$ we
know $u_n \to u$ and $ B(u_n) \to v$ in $W^{s,p}(I,F)$. Finally, since $B
\colon W^{s,p}(I,F) \to W^{s,p}(I,F)$
is demicontinuous, $B(u_n) \rightharpoonup B(u)$ in $W^{s,p}(I,F)$ and thus
$B(u) = v$.
\end{remark}
Our geometric condition is as follows:
\begin{itemize}
\item[$(A_5)$] [{\bf Monotonicity condition}] For each $ u \in {L^p(I,E)} $,
$$
    \alpha(A(u), (A+B)(u)) + \alpha( B(u), (A+B)(u)) \ge
\epsilon_0(L^p(I,E))
$$
\end{itemize}
\begin{remark} Condition $(A_5)$ is easier to verify that it might seem. For
instance, if $f = f(t)$ is constant, and then equation
(\ref{eqn:4}) reduces to a nonlinear generalization of the
Volterra equation, condition $(A_5)$ is immediately verified. In
this case all one needs is the reflexivity of $E$. Another common
way to verify condition $(A_5)$ is to assure the existence of a
cone $C$ in $L^p(I,E)$ with opening $2 - \epsilon_0(L^p(I,E))$
such that $\mathcal{I}m (A) \subseteq C$ and $ \mathcal{I}m (B)
\subseteq - C$. Indeed, let $ \zeta_A \in \mathcal{I}m (A)
\setminus \{ 0 \} \subseteq C$ and $\zeta_B \in \mathcal{I}m (B)
\setminus \{ 0 \} \subseteq (-C)$. There holds,
$$
\alpha(\zeta_A,\zeta_B) := \left \| \dfrac{\zeta_A}{\|\zeta_A\|} -
\dfrac{\zeta_B}{\|\zeta_B\|} \right \|  =
\left \| \dfrac{\zeta_A}{\|\zeta_A\|} +  \dfrac{\zeta_B}{\|\zeta_B\|}  - 2
\dfrac{\zeta_B}{\|\zeta_B\|} \right
\| \ge  2 - \alpha ( \zeta_A, -\zeta_B)  \ge \epsilon_0.
$$
To conclude we recall that for any three nonzero vectors,
$v_1,v_2,v_3$ in a normed vector space, we have, $ \alpha
(v_1,v_3) \le \alpha(v_1,v_2) + \alpha(v_2,v_3).$ Therefore
$$
\epsilon_0 \le \alpha(\zeta_A,\zeta_B) \le \alpha(\zeta_A,\zeta_A+ \zeta_B)
+ \alpha(\zeta_B, \zeta_A+\zeta_B ).
$$
\end{remark}
    We now can state the main result of this section.
\begin{theorem} \label{Existence} Assume ({$A_1$})--({$A_5$}).
Then there exists a $u \in L^p(0,T, E)$ solving the nonlinear integral
equation (\ref{eqn:4}).
\end{theorem}
    Before proving Theorem \ref{Existence} we need to develop a sequentially
weak continuity result for a class of operators acting on
vector-valued Lebesgue spaces. It is worthwhile to highlight that
such a result has it own importance in the theoretical point of
view. This is the content of what follows. The first result we
shall need is the fact that weak convergence in $L^\infty(E)$
implies a.e. weak convergence in $E$. This result was firstly
proven in \cite{K}. His proof though makes use of lifting theory
which brings an abstract flavor to it. The proof we shall present
here was taken from \cite{teixeira}, where some other consequences
of this fact is explored. Such a proof seems to give more insight
as to why this phenomenon should hold. Moreover, we also believe
the strategy presented here might be used in more general
situation to extract, almost everywhere (weak) convergence from
particular weakly convergent sequences, for instance minimizing
sequences in optimization problems. \\

\begin{theorem} \label{weakLinfty} Let ${E}$ be an arbitrary Banach space and $(\Omega, \mu)$
 be a Radon measure
 space. Let $u_n$ be a sequence in $L^\infty(\Omega, {E})$. Suppose $u_n \rightharpoonup u$
  in $L^\infty(\Omega,
{E})$. Then for $\mu$-almost every $x \in \Omega$,
$$
    u_n(x) \rightharpoonup u(x)  \textrm{ in } {E}.
$$
\end{theorem}
\begin{proof} Let $\varphi \in {E}^{*}$ be fixed. For each $ x \in \Omega $ and
$ 0 < r < \textrm{dist}(x ,
\partial \Omega)$, we define $\Phi_r^x \in \big [ L^\infty( \Omega, {E} ) \big ]^{*}$ to be
\begin{equation} \label{Phi}
    \Phi_r^{x} (f) := \dfrac{1}{\mu(B(x,r))} \int_{B(x,r)} \varphi(f(\xi)) d\mu(\xi).
\end{equation}
We verify that
$$
    \|\Phi_r^{x} \|_{ [ L^\infty( \Omega, {E} ) ]^{*} } :=
    \sup\limits_{ f \in L^\infty( \Omega, {E}
    )\setminus\{0\}} \dfrac{\Phi_r^{x}(f)}{\| f \|_{ L^\infty( {E} )}} \le
     \| \varphi \|_{{E}^{*}}.
$$
It follows therefore, from the Banach-Alaoglu Theorem, that, for a
fixed $x \in \Omega$, up to a subnet, we have
$$
    \Phi_r^{x} \stackrel{\star}{\rightharpoonup} \Phi^x \in \big [ L^\infty( \Omega, {E} )
    \big ]^{*}.
$$
Let us, hereafter, denote $v_n := u_n - u \in L^\infty( \Omega,
{E} )$. Let $\Omega_n$ be the Lebesgue set of $v_n$ provided by
Lebesgue's differentiation Theorem . We then set $ \Omega_0 =
\bigcap\limits_{n=1}^\infty\Omega_{n}$. In this way, $\Omega_0$
has total measure and for each $ x \in \Omega_0$ there holds
$$
    \Phi_r^{x}(v_n) \stackrel{r\to 0}{\longrightarrow} \varphi(v_n(x)) =
    \Phi^x (v_n(x)) \stackrel{n\to \infty}{\longrightarrow} 0.
$$
\end{proof}

In our next result we shall make use of Dunford's Theorem, which
we will state for convenience.
\begin{theorem}[Dunford]\label{D} Let $(\Omega, \Sigma, \mu)$ be a finite
measure space and ${X}$ be a
Banach space such that both ${X}$ and ${X}^{*}$ have the Radon-Nikod\'ym
property. A subset $K$ of $L^1(\Omega,
{X})$ is relatively weakly compact if and only if
\begin{enumerate}
\item $K$ is bounded, \item $K$ is uniformly integrable, and
\item for each $B \in \Sigma$, the set $\{ \int_B f
d\mu: f \in K \}$ is relatively weakly compact.
\end{enumerate}
\end{theorem}
\begin{lemma} \label{weak continuity}Let $p,~ q\ge 1$ and $I\colon
L^p(\Omega, {E}) \to L^\infty(\Omega,{E})$ be a
continuous linear map. Let $f\colon\Omega \times {E} \to {E}$ be a weak
Carath\'eodory map satisfying
$$
    \|f(x,u)\|_{{E}} \le A(x)\psi(\|u\|_{{E}}),
$$
where $A \in L^q(\Omega)$ and $\psi \in L^\infty_{\mathrm{loc}}(\Omega)$.
Then if either $q>1$ or $p=q=1$, the
map $\Psi := N_f \circ I \colon L^p(\Omega, {E}) \to L^q(\Omega, {E})$ is
sequentially weakly continuous.
\end{lemma}
\begin{proof} Let us suppose $q>1$. Let $u_n \rightharpoonup u$ in
$L^p(\Omega, {E})$. Since $\Psi$ is a bounded operator and
$L^q(\Omega, {E})$ is reflexive, up to a subsequence, $\Psi(u_n)
\rightharpoonup v \in L^q(\Omega, {E})$ for some $v \in
L^q(\Omega, {E})$. The idea is to show that actually $v =
\Psi(u)$. From Theorem \ref{weakLinfty}, we know $I(u_n)(x)
\rightharpoonup I(u)(x)$ in $E$ for $\mu$-a.e. $x\in \Omega$.
Since $f$ is a weak Carath\'eodory map, $\Psi(u_n)(x)
\rightharpoonup \Psi(u)(x)$ in $E$ for $\mu$-a.e. $x\in \Omega$ as
well. Now we shall conclude that $v = \Psi(u)$ $\mu$-a.e. To this
end we start by throwing away a set $\mathcal{A}_0$ of measure
zero such that
$$
F := \overline{ \textrm{span} }\Big [v \big ( \Omega \setminus \mathcal{A}_0
\big ) \cup \Psi(u) \big ( \Omega
\setminus \mathcal{A}_0 \big ) \Big ]
$$
is a separable and reflexive Banach space. The existence of such a
$\mathcal{A}_0$ is due to Pettis' Theorem.
Let now $\{\varphi_j\}$ be a dense sequence of continuous linear functionals
in $F$. By Ergorov's Theorem, for
each $\varphi_j$ fixed, there exists a negligible set $\mathcal{A}_j$, such
that $\varphi_j(v) =
\varphi_j(\Psi(u))$ in $\Omega \setminus \mathcal{A}_j$. Finally we define
$\mathcal{A} =
\bigcup\limits_{j=0}^\infty \mathcal{A}_j$. In this way $ \mu(\mathcal{A})
= 0$ and by the Hahn-Banach Theorem, $ v(x) = \Psi(u)(x)$ for all $x \in
\Omega \setminus \mathcal{A}$. \\
\indent Let us now study the case when $p=q=1$. For simplicity, we will
restrict ourselves to finite measure
spaces. We shall use Dunford's Theorem. Let $u_n \rightharpoonup u$ in
$L^1(\Omega, E)$. By the Eberlein-Smulian
Theorem the set $K = \{u, ~ u_n\}_{n=1}^\infty$ is weakly compact. Let us
show $\Psi(K)$ is relatively weakly
compact in $L^1(\Omega, E)$. Clearly $\Psi(K)$ is bounded, once
$\|\Psi(v)\|_{L^1(\Omega, E)} \le \|A\|_{L^1}
\cdot \psi(\|I\|\cdot \|v\|_{L^1(\Omega, E)})$. The last inequality also
shows $\Psi(K)$ is uniformly
integrable. Since $E$ is reflexive, we get item (3) of Durford's Theorem for
free. Hence, $\Psi(K)$ is
relatively weakly compact in $L^1(\Omega, E)$. Now we proceed as in the
previous case.
\end{proof}
We now have all the ingredients to prove Theorem \ref{Existence}.\\
{\it Proof of Theorem \ref{Existence}. } As done before, let
$g(t,x) := f(t,x) - f(t,0)$ and $B \colon L^p(I,E) \to L^p(I,E)$
be $B(u)(t) := g(t,u(t))$. We estimate
\begin{equation} \label{B}
\begin{array}{lll}
\|B(u)\|_{p} &=&  \left ( \displaystyle \int_0^T \| g(t,u(t)) \|_E^pdt
\right )^{1/p} \\
            &\le & \|u\|_{p}.
            \end{array}
\end{equation}
Let $\mathcal{K}\colon L^p(0,T,E) \to L^\infty(0,T,E)$ be the following map
$$
\mathcal{K}(u)(t) := \int_0^t k(t,\lambda)u(\lambda)d\lambda.
$$
We estimate
$$
\begin{array}{lll}
    \|\mathcal{K}(u)(t)\|_{E} &\le& \displaystyle\int_0^t |k(t,\lambda)|
\cdot \|u(\lambda)\|_{{E}}
     d\lambda \vspace{.2cm} \\
    &\le& \|k(t,\cdot)\|_q \cdot \|u\|_p\vspace{.2cm} \\
    &\le& C \|u\|_p.
\end{array}
$$
Above inequality says that
\begin{equation} \label{Estimate on K}
\|\mathcal{K}(u)\|_{\infty} \le C \|u\|_p.
\end{equation}
We remind  $A \colon L^p(I,E) \to L^p(I,E)$ was defined to be
$$
A (u)(t):= f(t,0) + \Phi \left (t, \int_0\sp t k(t,s)u(s)ds \right ).
$$
Showing the existence of a solution to problem (\ref{eqn:4}) is equivalent
to finding a fixed point to $A+B$.
Let us now estimate $\|A(u)\|_p$~:
\begin{equation} \label{Estimate on A}
\begin{array}{lll}
    \|A(u)\|_p &\le& \|f(\cdot,0)\|_p + \| \Phi(t,\mathcal{K}(u)(t))\|_p \\
    & \le & \|f(\cdot,0)\|_p + \|G\|_p\cdot \psi(C\cdot \|u\|_p).
\end{array}
\end{equation}
Let $M :=B_{L^p({E})}(\overline{R})$, be the ball in $L^p(0,T,{E})$ with
radius $\overline{R}$. It follows from
(\ref{Estimate on A}) that if $\|u\|_p \le \overline{R}$ then
$$
    \|A(u)\|_p \le \|f(\cdot,0)\|_p + \|G\|_p\cdot \psi(C\cdot \overline{R})
\le \overline{R}.
$$
It implies that $A$ maps $M$ into itself. Estimate (\ref{B}) shows
that $B$ also maps $M$ into itself. Moreover Lemma \ref{weak
continuity} says $A\colon M \to M$ is sequentially weakly
continuous. We now make use of the monotonicity condition to
estimate the size of $ (A+B)(M) $. From Lemma \ref{Triang},
$$
\| A(u) + B(u) \| \le 2 \overline{R} \Big [ 1 - \big [ \alpha(A(u),A(u) +
B(u)) + \alpha(A(u) + B(u),B(u) ) \big
] \Big ] \le \overline{R}.
$$
Keeping in mind that $M$ is a compact subset of a locally convex Hausdorff
space, and $A+B$ is a continuous map
from $M$ into itself, we guarantee, by the Schauder-Tychonoff Theorem, the
existence of a fixed point to $A+B$.
\hfill $\Box$
\begin{remark} An important case of equation (\ref{Existence}) is
when $f \in L^\infty(I,E)$, $\| k(t,u)\| \le C$ and $\|
\Phi(t,u)\| \le G\psi(\|u\|)$. In this case, as we have observed
above, we do not need the geometric condition $(A_5)$. By using a
similar idea of section 6, we can build the locally convex
topology in $L^\infty(I,E)$ given by $\bigcup\limits_{n\ge 2}
\mathcal{T}^n$, where $\mathcal{T}^n := \{ \rho\colon
L^\infty(I,E) \to \Bbb{R}_+ : \rho = |f| \text{ and } f \in
[L^n(I,E)]^{*} \}$. In this topology the ball of $L^\infty(I,E)$
is compact. Furthermore we can also verify the continuity of the
operators $A$ and $B$ with respect to this new topology. Thus, in
this particular case, we can find a bounded solution of equation
(\ref{Existence}). This is as much as one should expect due to the
nonlinearities involved in the equation.
\end{remark}

\end{document}